\documentclass{amsart}

\usepackage{amssymb,amsmath}
\usepackage[neveradjust]{paralist}
\usepackage{epic,eepic,color}
\usepackage[all]{xypic}
\xyoption{dvips}
\usepackage{url}
\usepackage[colorlinks,plainpages,backref]{hyperref}
\usepackage{dsfont}
\renewcommand{\mathbb}{\mathds}
\usepackage{mathrsfs}
\DeclareMathAlphabet{\mathsc}{U}{rsfs}{m}{n}
\renewcommand{\mathcal}{\mathsc}

\theoremstyle{definition}
\newtheorem{ntn}{Notation}[section]
\newtheorem{dfn}[ntn]{Definition}
\newtheorem*{ack}{Acknowledgment}
\theoremstyle{plain}
\newtheorem{lem}[ntn]{Lemma}
\newtheorem{prp}[ntn]{Proposition}
\newtheorem{thm}[ntn]{Theorem}
\newtheorem{cor}[ntn]{Corollary}

\theoremstyle{remark}

\newtheorem{exa}[ntn]{Example}

\numberwithin{equation}{section}

\renewcommand{\aa}{\mathbf{a}}
\newcommand{\A}{\mathcal{A}}
\renewcommand{\AA}{\mathbb{A}}

\newcommand{\E}{\mathcal{E}}
\newcommand{\F}{\mathcal{F}}

\newcommand{\M}{\mathcal{M}}
\newcommand{\C}{\mathbb{C}}
\newcommand{\Q}{\mathbb{Q}}
\newcommand{\PP}{\mathbb{P}}
\renewcommand{\SS}{\mathcal{S}}
\newcommand{\pp}{\mathfrak{p}}
\newcommand{\p}{\partial}
\newcommand{\ideal}[1]{{\left\langle#1\right\rangle}}
\newcommand{\into}{\hookrightarrow}
\newcommand{\oma}{\omega_{\aa}}
\newcommand{\onto}{\twoheadrightarrow}
\renewcommand{\O}{\mathcal{O}}
\newcommand{\xymat}{\SelectTips{cm}{}\xymatrix}

\newcommand{\ZZ}{\mathbb{Z}}
\newcommand{\Z}{\mathcal{Z}}

\newcommand{\set}[1]{\left\{#1\right\}}
\newcommand{\abs}[1]{\left|#1\right|}

\DeclareMathOperator{\Ext}{Ext}
\DeclareMathOperator{\SExt}{{\mathcal Ext}}
\DeclareMathOperator{\Hom}{Hom}
\DeclareMathOperator{\pd}{pd}
\DeclareMathOperator{\Proj}{Proj}

\DeclareMathOperator{\Sing}{Sing}
\DeclareMathOperator{\Spec}{Spec}
\DeclareMathOperator{\Sym}{Sym}
\DeclareMathOperator{\td}{td}
\DeclareMathOperator{\length}{length}
\DeclareMathOperator{\codim}{codim}
\DeclareMathOperator{\ch}{ch}

\begin{document}

\title[Complexes of logarithmic forms along arrangements]
{Complexes, duality and Chern classes of logarithmic forms along hyperplane arrangements}

\author{Graham Denham}
\address{
G.~Denham\\
Department of Mathematics\\
Middlesex College\\
University of Western Ontario\\
London, ON N6A 5B7\\
Canada
}
\email{gdenham@uwo.ca}
\thanks{GD was partially supported by a grant from NSERC of Canada.}

\author{Mathias Schulze}
\address{
M.~Schulze\\
Department of Mathematics\\
Oklahoma State University\\
Stillwater, OK 74078\\
United States}
\email{mschulze@math.okstate.edu}
\thanks{}

\date{\today}

\begin{abstract}
We describe dualities and complexes of logarithmic forms and differentials for central affine and corresponding projective arrangements.
We generalize the Borel--Serre formula from vector bundles to sheaves on $\PP^d$ with locally free resolutions of length one.
Combining these results we present a generalization of a formula due to Musta{\c{t}}{\u{a}} and Schenck, relating the Poincar\'e polynomial of an arrangement in $\PP^3$ (or a locally tame arrangement in $\PP^d$ with zero-dimensional non-free locus) to the total Chern polynomial of its sheaf of logarithmic $1$-forms.
\end{abstract}



\maketitle
\tableofcontents

\section{Introduction}

The study of logarithmic differentials and differential forms in the context of hyperplane arrangement singularities now has a thirty-year history that begins with Saito~\cite{Sai80}.  For a comprehensive 
survey, we refer to the forthcoming book \cite{CDFSSTY}.  This
paper collects together some results about complexes of logarithmic forms: some known, others folklore, and some new.  In particular, we attempt to clarify the relationship between forms on the complements of central, non-central and projective arrangements, respectively.  

In order to be more precise, let $\Omega^p(\A)$ denote the module of logarithmic $p$-forms (defined in \eqref{eq:defOmega}), for a central arrangement $\A$ of rank $\ell$. 
Following \cite{Dol82,Dol07,Dim92}, we consider the submodule $\Omega^p_0(\A)$ (Definition~\ref{25}) of forms which vanish when contracted along the Euler differential. 
We show that the coherent sheaf on this submodule coincides with the usual sheaf of logarithmic forms on the projectivization of $\A$, Definition~\ref{dfn:Omega-proj}. 
We also consider a {\em quotient} of $\Omega^p(\A)$ which we call the module of {\em relative} logarithmic forms, denoted $\Omega_\sigma^p(\A)$ and defined in Definition~\ref{dfn:omegasigma}. 
This was introduced in slightly different terms by Terao and Yuzvinsky~\cite{TY95}, and discussed more generally in \cite{GMS09}. 
We note that this module is isomorphic to $\Omega_0^p(\A)$, which gives a noncanonical splitting of the inclusion $\Omega_0^p(\A)\hookrightarrow\Omega^p(\A)$. 
We note that the choice involved amounts to choosing an affine chart. 
In this way, we understand logarithmic forms for affine, non-central arrangements in terms of their cone.

The modules of logarithmic forms are reflexive, which was observed first in \cite[(1.7)~Cor.]{Sai80} for $\Omega^1(\A)$.  
The exterior product gives a map $j\colon \bigwedge^p\Omega^1(\A)\to\Omega^p(\A)$, which we see is an isomorphism exactly when $\bigwedge^p\Omega^1(\A)$ is also reflexive (Proposition~\ref{prp:doubledual}.) 
We note that, if $\A$ is tame, (Definition~\ref{df:tame}), then $j$ is an isomorphism for values of $p$ less than the codimension of the non-free locus of $\A$ (Proposition~\ref{prp:Lebelt}). 
In the case where $\A$ is free, this is part of Saito's criterion for freeness, from \cite[(1.8)~Thm.]{Sai80}; in the case where $\A$ is locally free, it was noted by Musta{\c{t}}{\u{a}} and Schenck~\cite{MuSc01}. 
The same is true for the variations above.

In \S\ref{sec:dual}, we examine the duals of the modules of logarithmic forms, which are modules (or
sheaves) of multilinear logarithmic differentials.  In the dual setting, the natural construction is a quotient of the module of logarithmic forms $D_p(\A)$, denoted by $D^0_p(\A)$ and defined in Definition~\eqref{dfn:defD0}. 
Our work with forms allows us to replace the quotient with a submodule, $D^\sigma_p(\A)$, again by choosing a chart. 
The modules of forms are also self-dual, which gives some useful symmetry.  We note some equivalent formulations of the homological notions of free and tame arrangements.

Multiplication by a degree-$0$ logarithmic $1$-form $\omega_\lambda$
gives a cochain complex $(\Omega^\bullet(\A),\omega_\lambda)$,
as well as for the projective constructions.  Orlik and Terao~\cite{OT95} show that this complex is exact, and  $(\Omega^\bullet_\sigma(\A),\omega_\lambda)$ has a single non-zero cohomology group, both for suitably ``generic'' choices of $\lambda$, which we discuss below. 
Their main application is to show that, under the same hypotheses on $\lambda$, the function
\begin{equation}\label{eq:master}
\Phi_\lambda=\prod_{H\in\A}\alpha_H^{\lambda_H}
\end{equation}
has isolated, nondegenerate critical points.  (Here, $\A$ is a non-central arrangement defined by affine-linear forms $\set{\alpha_H:H\in\A}$.)  Recent generalizations of this result in \cite{CDFV09} make use of a parameterized version of this complex, which we discuss in \S\ref{ss:parametric}.  We clarify the relationship between the various complexes and, in doing so, improve slightly on a main result of \cite{CDFV09}.

As further applications, we give some formulas for the Chern classes of the sheaf of logarithmic $1$-forms in \S\ref{sec:chern}. 
Musta{\c{t}}{\u{a}} and Schenck~\cite{MuSc01} showed that, for locally free arrangements, the total Chern polynomial is given by the arrangement's characteristic polynomial. 
We see in Proposition~\ref{prp:slice-Chern} that the coefficients of the two polynomials always agree up to degree $k$, for any arrangement, if its non-free locus has codimension $>k$. 
For certain projective arrangements with zero-dimensional non-free locus (including all arrangements of rank $\leq4$), we compute the total Chern polynomial in Theorem~\ref{th:Chern}. 
Our expression is combinatorial if the non-free locus consists of generic closed subarrangements; however, we exhibit two arrangements with the same matroid and different Chern polynomials in Example~\ref{ex:ziegler}.

\section{Log forms}

\subsection{Normal and reflexive sheaves}

Let $\A$ be a central (simple) arrangement of $n$ hyperplanes in an $\ell$-dimensional $\C$-vector space $V=\AA^\ell$.

We denote by $L(\A)$ the intersection lattice of $\A$, and by $L_c(\A)\subseteq L(\A)$ the sublattice of codimension-$c$ flats.
For $X\in L(\A)$, let $\A_X=\{H\in\A\mid X\subseteq H\}$ be the localization of $\A$ at $X$; for $H\in\A$ let $\A^H=\{H'\cap H\mid H'\in\A\setminus\{H\}\}$ be the restriction of $\A$ to $H$.

We abbreviate $\PP^d=\PP V$, $d:=\ell-1$, and denote by $\PP\A:=\pi(\A)$ the corresponding projective arrangement, where 
\begin{equation}\label{32}
\pi\colon V\setminus\set{0}\to\PP V
\end{equation}
is the canonical projection. 
For $H\in\A$ we denote by $\alpha_H\in V^*$ its defining equation.
We can assume that for some $\{H_1,\dots,H_\ell\}\subset\A$, $x_i:=\alpha_{H_i}$ are coordinates on $V$.
Then $f=\prod_{H\in\A}\alpha_H$ is the (reduced) defining equation of $\A$ in the coordinate ring $R=\Sym(V^*)$ of $V$.

We shall readily identify $R$-modules with the associated sheaves on $V$, and denote $-^\vee=\Hom_R(-,R)$.
We identify $\A$ with $\Spec(R/f)\subset V$ and denote by $i\colon U:=V\backslash\Sing\A\into V$ the inclusion of the complement of the singular locus of $\A$ in $V$.
Following Barth~\cite{Ba77}, we say an $R$-module $M$ $\A$-normal if $M\to i_*i^*M$ is an isomorphism.
Note that $R$ itself is normal.
The following easy facts are particularly useful for our purpose.

\begin{lem}\label{31}\
\begin{enumerate}
\item\label{31a} Any reflexive module has depth $\geq2$.
\item\label{31b} Two $\A$-normal modules are dual if they are dual on $U$.
\item\label{31c} A torsion-free $\A$-normal module that is reflexive on $U$ is reflexive.
\item\label{31d} Moreover, $\A$-normality is a consequence of reflexivity.
\end{enumerate}
\end{lem}

\begin{proof}
The non-obvious statements follow from \cite[Props.~1.1, 1.6]{Har80}.
\end{proof}

\subsection{Relative logarithmic forms}

The modules $\Omega^p_V$ form a graded $R$-module $\Omega^\bullet_V$ closed under exterior product such that the natural map $R\to\Omega^1_V$ has degree $0$.
In particular, the localization $\Omega^\bullet_{V,f}$ has a natural $R$-grading.
The module of {\em logarithmic differential $p$-forms along $\A$} is the graded $R$-module
\begin{equation}\label{eq:defOmega}
\Omega^p(\A)=\left\{\omega\in\frac1f\Omega^p_V\mid\forall H\in\A\colon\frac{d\alpha_H}{\alpha_H}\wedge\omega\in
\frac1f\Omega_V^{p+1}\right\}\subset\Omega^\bullet_{V,f}.
\end{equation}
It is easy to check that $\Omega^\bullet(\A)$ is closed under exterior
product.

As the modules $\frac1f\Omega_V^p$ and $\frac1f\Omega_V^{p+1}$ are free, $\Omega^p(\A)$ is $\A$-normal, and 
\begin{equation}\label{13}
\Omega^p(\A)_\pp=\bigwedge^p\Omega^1(\A)_\pp
\end{equation}
for primes $\pp$ supported on the non-singular locus $U$. 
Both sides are free over $R_\pp$, so $\Omega^p(\A)$ is reflexive by Lemma~\ref{31}.\eqref{31c}.

The exterior product gives a map $j_p\colon \bigwedge^p\Omega^1(\A)\to\Omega^p(\A)$, which is easily seen to be a monomorphism. 
Under some hypotheses, $j$ is an isomorphism, as we will see in Proposition~\ref{prp:Lebelt}. 
However, exterior powers of a reflexive module need not be reflexive, and in general one has:

\begin{prp}\label{prp:doubledual}
For any $\A$, we have $(\bigwedge^p\Omega^1(\A))^{\vee\vee}\cong\Omega^p(\A)$ for $0\leq p\leq \ell$.
\end{prp}

\begin{proof}
Let $E^p=E^p(\A)$ denote the cokernel of $j_p$. 
By \eqref{13}, the module $E^p$ is supported on $\Sing\A$, which has codimension $\geq2$.
Therefore $\Ext^i_R(E^p,R)=0$ for $i=0,1$. 
It follows that
\begin{equation}\label{eq:wedgedual}
j_p^\vee\colon \Omega^p(\A)^\vee\to(\bigwedge^p\Omega^1(\A))^\vee
\end{equation}
is an isomorphism. 
We apply $-^\vee$ again and note $\Omega^p(\A)$ is reflexive.
\end{proof}

By Lemma~\ref{31}.\eqref{31a}, $\pd_R\Omega^p(\A)\leq\ell-2$ for all $p$. 
The following definition first appeared in \cite{OT95b}:

\begin{dfn}\label{df:tame}
An arrangement $\A$ is {\em tame} if $\pd_R\Omega^p(\A)\leq p$, for $0\leq p\leq\ell$.
\end{dfn}

Since $\Omega^0(\A)=R$ is free, the condition is vacuous except for $1\leq p\leq\ell-3$. 
Accordingly, all arrangements of rank $\leq3$ are tame.

The graded $R$-dual of $\Omega^1(\A)$ is the graded $R$-module of {\em logarithmic differentials along $\A$} 
\[
D(\A)=D_1(\A)=\{\delta\in D_V\mid\forall H\in\A\colon\delta(\alpha_H)\in\alpha_HR\}
\]
where $D_V={\rm Der}_\C(R,R)$ is the module of polynomial vector fields on $V$.

Note that the standard Euler differential 
$
\chi=\sum_{i=1}^\ell x_i\p_i
$
is a member of $D(\A)$.
Recall from \cite[Prop.~4.86]{OTbook} that contraction $\iota_\delta$ with $\delta\in D(\A)$ defines a graded map
\[
\iota_\delta\colon\Omega^p(\A)\to\Omega^{p-1}(\A).
\]
By $\A$-normality of $\Omega^\bullet(\A)$, it suffices to check this on $U$, where it is clear.

\begin{dfn}\label{25}
We call $\Omega^p_0(\A):=\ker\iota_\chi\subset\Omega^p(\A)$ be the module of {\em relative logarithmic differential $p$-forms along $\A$}.
\end{dfn}

Clearly, $\Omega^\bullet_0(\A)$ is closed under exterior product, and locally free on $U$.
As above, freeness of $\frac1f\Omega_V^p$ and $\frac1f\Omega_V^{p-1}$ implies that $\Omega^p_0(\A)$ is $\A$-normal and torsion free. 
Then by Lemma~\ref{31}.\eqref{31c} and the same argument as Proposition~\ref{prp:doubledual},

\begin{prp}\label{prp:Omref}
$\Omega^\bullet_0(\A)$ is closed under exterior product and reflexive.
Moreover, $\Omega^\bullet_0(\A)\cong(\bigwedge^\bullet\Omega^1_0(\A))^{\vee\vee}$.
\end{prp}

We note also the following simple fact (cf.~\cite[Exa.~4.122]{OTbook}).

\begin{prp}
$\A\mapsto\Omega^p_0(\A)$ is a local functor.
That is, for any graded prime ideal $\pp$ of $R$, $\Omega_0^p(\A)_\pp\cong\Omega_0^p(\A_{X(\pp)})_\pp$, where $X(\pp)$ denotes the set-theoretically smallest subspace in the intersection lattice $L(\A)$ containing the zeroes of $\pp$.
\end{prp}

The term ``relative'' in Definition~\ref{25} does not refer to a specific map here.
It turns out that the so defined differential forms are ``relative'' to many maps simultaneously.
This can be seen as follows:
As in \cite{CDFV09}, we denote 
\begin{equation}\label{33}
\omega_\aa=d\log(\alpha^\aa)=\sum_{H\in\A}a_H\frac{d\alpha_H}{\alpha_H}
\end{equation}
for independent variables $\aa=(a_H)_{H\in\A}\in(\C^\A)^\vee$.
We write 
\begin{equation}\label{28}
|\aa|=\sum_{H\in\A}a_H=\iota_\chi(\omega_\aa).
\end{equation}
Now for any $\sigma\in\C^\A$ with $|\sigma|\ne0$, $(\Omega^\bullet(\A),\omega_\sigma)$ is split exact (see \cite[Prop.~4.86]{OTbook}), since
\[
\iota_\chi(\omega_\aa\wedge\omega)=|\aa|\omega-\omega_\aa\wedge\iota_\chi(\omega).
\]

\begin{dfn}\label{dfn:omegasigma}
We define module of {\em logarithmic differential $p$-forms relative to $\alpha^\sigma$} to be the graded $R$-module
\[
\Omega^p_\sigma(\A)=\Omega^p(\A)/\omega_\sigma\wedge\Omega^{p-1}.
\]
\end{dfn}

For any $\lambda\in\C^\A$ with $|\lambda|=0$, we can identify the complexes 
\begin{equation}\label{3}
(\Omega^\bullet_\sigma(\A),\omega_\lambda)=(\Omega^\bullet_0(\A),\omega_\lambda)
\end{equation}
via the map $[\omega]\mapsto\omega-\omega_\sigma\wedge\iota_\chi(\omega)$.
In particular, we can identify
\begin{equation}\label{1}
\Omega^\bullet(\A)=\omega_\sigma\wedge\Omega^{\bullet-1}_0(\A)\oplus\Omega^\bullet_0(\A)\cong\Omega^{\bullet-1}_0(\A)\oplus\Omega^\bullet_0(\A).
\end{equation}
Then the notion of tameness, Definition~\ref{df:tame}, is unaffected by working with relative logarithmic forms:

\begin{prp}\label{prp:tame0}
$\A$ is tame if and only if $\pd_R\Omega^p_0(\A)\leq p$ for $1\leq p\leq \ell-3$.
\end{prp}

\begin{proof}
If $\pd_R\Omega^p_0(\A)\leq p$ for all $p$, then the same is true for $\Omega^p(\A)$ by \eqref{1}. 
Conversely, if $\A$ is tame, argue by induction on $p$ using \eqref{1} that $\Ext^q_R(\Omega_0^p(\A),R)=0$ for all $q>p$.  
\end{proof}

The following generalizes a result of Musta{\c{t}}{\u{a}}--Schenck~\cite[Thm.~5.3]{MuSc01}.

\begin{prp}\label{prp:Lebelt}
If $\A$ is a tame arrangement whose non-free locus has codimension $>k$, then $\bigwedge^p\Omega^1_0(\A)=\Omega^p_0(\A)$, and hence also $\bigwedge^p\Omega^1(\A)=\Omega^p(\A)$, for $p<k$.
\end{prp}

\begin{proof}
By the argument from the proof of \cite[Lem.~5.2]{MuSc01}, the hypotheses imply that $\Omega^1_0(\A)$ is a $(k-1)$-syzygy.
Then, by \cite[Satz~3.1]{Leb77}, $\pd_R\bigwedge^p\Omega^1_0(\A)=p$ for $p\le k$.
Applying $\Ext_R^q(-,R)$ to the exact sequence
\[
0\to\bigwedge^p\Omega^1_0(\A)\to\Omega^p_0(\A)\to E^p_0\to 0
\]
shows that $\Ext_R^q(E^p_0,R)=0$ for $q>p+1$ and hence for $q>k>p$.
But by assumption the $E^p_0$ are supported in codimension $>k$ only,
so $\Ext^q_R(E^p_0,R)=0$ for $q\leq k$ as well, so $E^p_0=0$ for
all $p<k$.
\end{proof}

\subsection{Projective logarithmic forms}

Geometrically, $\Omega_0^\bullet(\A)$ can be considered as an affine version of the sheaf $\Omega^\bullet(\PP\A)$ defined following the original definition of K.~Saito~\cite{Sai80}, as we will see in Proposition~\ref{prp:proj}.
Recall that $\PP\A$ is the image of $\A$ under the natural projection \eqref{32}.

\begin{dfn}\label{dfn:Omega-proj}
Consider $\PP\A$ as a principal divisor $D=(f)$ on $\PP V$.
Then we call
\[
\Omega^p(\PP\A):=\Omega_{\PP V}^p(D)\cap d^{-1}\Omega_{\PP V}^{p+1}(D).
\]
is the sheaf of {\em logarithmic differential $p$-forms along $\PP\A$}.
\end{dfn}

In other words, $\Omega^\bullet(\PP\A)$ is the sheaf of rational differential forms $\omega$ on $\PP V$ for which both $\omega$ itself and $d\omega$ have at most a simple pole along $\PP\A$.
Note that in any chart $\{x_i=1\}$, $\Omega^\bullet(\PP\A)$ restricts to $\Omega^\bullet(\A^{\{x_i=1\}})$.

In order to see the claimed relation with $\Omega_0^\bullet(\A)$, first fix a chart with index $i\in\{1,\dots,\ell\}$.
Then $\Omega^\bullet_{V,f}/\omega_{e_i}\wedge\Omega^\bullet_{V,f}$ is the module of differential forms with poles along $\A\setminus\{H_i\}$ relative to the map $x_i$ on $\{x_i\ne0\}$.
Using the definition of $\Omega^\bullet(\A)$, one checks that the composition of canonical maps
\[
\Omega^\bullet(\A)_{x_i}\into\Omega^\bullet_{V,f}\onto\Omega^\bullet_{V,f}/\omega_{e_i}\wedge\Omega^\bullet_{V,f}
\]
factors through an inclusion
\[
\Omega_{e_i}^\bullet(\A)_{x_i}\into\Omega^\bullet_{V,f}/\omega_{e_i}\wedge\Omega^\bullet_{V,f},
\]
and it follows that we can consider
\begin{gather}\label{21}
\Omega_{e_i}^\bullet(\A)\otimes_RR/\ideal{x_i-1}=\Omega^\bullet(\A^{\{x_i=1\}}),\\
\nonumber\Omega_{e_i}^\bullet(\A)_{x_i}=\Omega^\bullet(\A^{\{x_i=1\}})\otimes_\C\C[x_i^{\pm1}].
\end{gather}
This immediately implies (see~\cite{TY95})

\begin{prp}\label{22}
This correspondence \eqref{21} combined with \eqref{3} identifies $(\Omega^\bullet_0(\A)_{x_i},\omega_\lambda)$ and $(\Omega^\bullet(\A^{\{x_i=1\}}),\omega_{\lambda_{\hat i}})$ where $\lambda_{\hat i}$ is obtained by deleting the $i$th component from $\lambda$. 
\end{prp}

In order to understand the global relation of $\Omega_0^\bullet(\A)$ and $\Omega^\bullet(\PP\A)$, we follow that approach in the proof of \cite[Thm.~8.4]{Har77}.
We obtain a logarithmic version of the well-known description of logarithmic differential forms on projective space as the kernel of contraction with the Euler differential, cf.~\cite[Ch.~6, \S1]{Dim92} and \cite[\S2.1]{Dol82}. 
For $1$-forms along a generic $\A$, this result can be found in \cite[p.~702-703]{MuSc01}.

\begin{prp}\label{prp:proj}
As sheaves on $\PP V$, 
\[
\widetilde{\Omega_0^p(\A)}=\Omega^p(\PP\A).
\]
\end{prp}

\begin{proof}
Analogous to the logarithmic version in Definition~\ref{25}, define
\begin{equation}\label{26}
\Omega^p_{V,0}:=\ker(\iota_\chi\colon\Omega_V^p\to\Omega^{p-1}_V).
\end{equation}
Consider the sequence
\[
\xymat{
0\ar[r]&\Omega_{V,0}^1\ar[r]&\Omega_V^1\ar[r]^-{\iota_\chi}&\Omega_V^0\ar[r]&0,
}
\]
which is exact away from the origin.
In the proof of \cite[Thm.~8.4]{Har77}, it is denoted $0\to M\to E\to S$.
Applying $-^\vee$, $\bigwedge^p$, and then $-^\vee$ again, yields a sequence
\begin{equation}\label{23}
\xymat{
0\to\bigwedge^p\Omega_{V,0}^1\ar[r]&\Omega_V^p\ar[r]^-{\iota_\chi}&
\Omega^{p-1}_V,
}
\end{equation}
which is exact away from the origin.
By \cite[Thm.~8.4]{Har77}, $\widetilde\Omega_{V,0}^1=\Omega^1_{\PP V}$, and hence 
\begin{equation}\label{24}
\widetilde{\bigwedge^p\Omega_{V,0}^1}=\Omega^p_{\PP V}.
\end{equation}
It follows from \eqref{26}, \eqref{23}, and \eqref{24} that $\widetilde\Omega^p_0=\Omega^p_{\PP V}$, and then
\[
\widetilde{\Omega^p_0(*\A)}=\Omega^p_{\PP V}(*\PP\A).
\]
Comparison of the subsheaves $\widetilde{\Omega^p_0(\A)}\subset\widetilde{\Omega^p_0(*\A)}$ and $\Omega^p(\PP\A)\subset\Omega^p_{\PP V}(*\PP\A)$ can now be done in charts.
So the claim follows from Proposition~\ref{22}.
\end{proof}

\begin{dfn}
We call $\PP\A$ \emph{locally tame} if $\Omega^p(\PP\A)$ has a locally free resolution of length $p$, for $0\leq p\leq\ell-1$.
\end{dfn}

All arrangements in $\PP^3$ are locally tame by reflexivity of $\Omega^1(\PP\A)$ and \cite[Prop.~1.3]{Har80}.
Proposition~\ref{prp:Lebelt} applied in charts gives the following:

\begin{prp}\label{prp:Lebelt-proj}
If $\PP\A$ is a locally tame arrangement whose non-free locus has codimension $>k$, then $\bigwedge^p\Omega^1(\PP\A)=\Omega^p(\PP\A)$ for $p<k$.
\end{prp}

Recall that $\A$ is called locally free if $\Omega^1(\A_X)$ is free for all $X\in L_{<\ell}(\A)$.

\begin{dfn}
We call $\PP\A$ locally free if the sheaf $\Omega^1(\PP\A)$ is a vector bundle.
\end{dfn}

From Proposition~\ref{prp:proj}, using that $\A\mapsto\Omega^1(\A)$ is a local functor, we deduce the following equivalence.

\begin{lem}\label{lem:locfree}
$\PP\A$ is locally free if and only if $\A$ is locally free.
More precisely, $\Omega^1(\PP\A)$ is free on an open set $U\subseteq\PP V$ if and only if, whenever $X\subseteq U$ for some $X\in L(\PP\A)$, the closed subarrangement $\A_X$ is free. 
\end{lem}

\section{Dualities}\label{sec:dual}

\subsection{Duality with relative log differentials}

We define the module of {\em logarithmic differentials relative to $\alpha^\sigma$} to be the graded $R$-module
\begin{equation}\label{7}
D^\sigma(\A):=\{\delta\in D(\A)\mid \delta(\alpha^\sigma)=0\}.
\end{equation}
Then clearly
\begin{equation}\label{11}
D(\A)=R\chi\oplus D^\sigma(\A).
\end{equation}
Applying $-^\vee$ to the exact sequence
\[
\xymat{
0\ar[r]&D^\sigma(\A)\ar[r]&D(\A)\ar[r]^-\phi&R\ar[r]&0,
}
\]
where $\phi(\delta)=\frac{\delta(\alpha^\sigma)}{\alpha^\sigma}=\langle\omega_\sigma,\delta\rangle$, shows that $\phi^\vee=\omega_\sigma$ hence

\begin{lem}
$D^\sigma(\A)^\vee=\Omega^1_\sigma(\A)$.
\end{lem}

The higher $R$-modules of logarithmic differentials are defined by
\begin{equation}\label{def-Dp}
D_p(\A)=\left\{\delta\in\bigwedge^pD_V\mid\forall H\in\A,g_2,\dots,g_p\in R\colon\delta(\alpha_H,g_2,\dots,g_p)\in\alpha_HR\right\}.
\end{equation}
As in \eqref{13}, we have 
\begin{equation}\label{12}
D_p(\A)_\pp=\bigwedge^p D_1(\A)_\pp
\end{equation}
for primes $\pp$ supported on $U$
Generalizing \eqref{7}, we introduce higher relative logarithmic differentials.

\begin{dfn}
We define the module of {\em logarithmic $p$-differentials relative to $\alpha^\sigma$} to be the graded $R$-module
\[
D_p^\sigma(\A)=\left\{\delta\in D_p(\A)\mid\forall g_2,\dots,g_p\in R\colon\delta(\alpha^\sigma,g_2,\dots,g_p)=0\right\}.
\]
\end{dfn}

Both $D_\bullet(\A)$ and $D_\bullet^\sigma(\A)$ are clearly closed under exterior product.
As in the case of forms, both are $\A$-normal and torsion free, and using Lemma~\ref{31} we obtain

\begin{prp}
$D_\bullet^\sigma(\A)$ is closed under exterior product and reflexive.
\end{prp}

By \cite[Prop.~2.2]{MuSc01}, there is a non-degenerate pairing
\begin{equation}\label{8}
\ideal{-,-}\colon
\Omega^p(\A)\times D_p(\A)\to R.
\end{equation}
By $\A$-normality of the three modules involved, it is sufficient to check this on $U$ where it is clear.

Under duality such as in \eqref{15}, subspaces of one factor correspond to quotient spaces of the other factor.
Therefore we need also a quotient representation of $D_p^\sigma(\A)$ independent of $\sigma$.

\begin{dfn}\label{dfn:defD0}
We define the module of {\em relative logarithmic $p$-differentials along $\A$} to be the graded $R$-module
\[
D_p^0(\A):=D_p(\A)/\chi\wedge D_{p-1}(\A).
\]
\end{dfn}

As we can check on $U$, $(D_\bullet(\A),\chi)$ is exact and is splits by $\iota_{\omega_\sigma}$.
We can thus, as in \eqref{3} and \eqref{1}, identify
\[
D_p^\sigma(\A)=D_p^0(\A),\quad D_p(\A)=\chi\wedge D_{p-1}^0(\A)\oplus D_p^0(\A)\cong D_{p-1}^0(\A)\oplus D_p^0(\A).
\]

\begin{prp}\label{prp:dual}
The pairing \eqref{8} induces a non-degenerate pairing
\begin{equation}\label{10}
\Omega^p_0(\A)\times D_p^0(\A)=\Omega^p_\sigma(\A)\times D_p^\sigma(\A)\to R.
\end{equation}
\end{prp}

\begin{proof}
For well-definedness of \eqref{10}, we need to show that $\big(\omega_\sigma\wedge\Omega^{p-1}(\A)\big)\times D_p^\sigma(\A)$ and $\Omega^p_0(\A)\times\big(\chi\wedge D_{p-1}^0(\A)\big)$ are mapped to zero by \eqref{8}, using the sub and quotient representations the two factors.
This can be checked on $U$, and we show only the first statement.
By \eqref{13} and \eqref{12}, locally at points in $U$, $\Omega^{p-1}(\A)$ and $D_p(\A)$ are generated by $\omega=\frac{d\alpha_{H_2}}{\alpha_{H_2}}\wedge\dots\wedge\frac{d\alpha_{H_p}}{\alpha_{H_p}}$ and $\delta=\delta_1\wedge\dots,\wedge\delta_p$, $\delta_i\in D(\A)$, respectively.
Using \eqref{33},
\[
\ideal{\omega_\sigma\wedge\omega,\delta}=\frac{\det(\delta_i(\alpha^\sigma)\vert\delta_i(\alpha_{H_j}))}{\alpha^\sigma\alpha_{H_2}\cdots\alpha_{H_p}}=\frac{\delta(\alpha^\sigma,\alpha_{H_1},\dots,\alpha_{H_p})}{\alpha^\sigma\alpha_{H_2}\cdots\alpha_{H_p}},
\]
and well-definedness follows.

Now we need to verify conversely that any $\omega\in\Omega^p(\A)$ with $\ideal{\omega,D_p^\sigma(\A)}=0$ must be in $\omega_\sigma\wedge\Omega^{p-1}(\A)$, and that 
any $\delta\in D_p(\A)$ with $\ideal{\Omega^p_0(\A),\delta}=0$ must be in $\chi\wedge D_{p-1}(\A)$.
Again we show only the first statement, and we can restrict ourselves to local considerations on $U$:
By exactness of $(\Omega^\bullet(\A),\omega_\sigma)$, it suffices to show that $\omega_\sigma\wedge\omega=0$.
From \eqref{11} and \eqref{12} we derive that
\[
D_{p+1}(\A)=(R\chi\oplus D^\sigma(\A))\wedge\bigwedge^p D^\sigma(\A).
\]
As $\ideal{\omega_\sigma\wedge\omega,\bigwedge^{p+1}D^\sigma(\A)}=0$ and, for $\delta\in D_p^\sigma(\A)$,
\[
\ideal{\omega_\sigma\wedge\omega,\chi\wedge\delta}=|\sigma|\ideal{\omega,\delta}=0
\]
by hypothesis, we find that $\ideal{\omega_\sigma\wedge\omega,D_{p+1}(\A)}=0$.
Using \eqref{12} and local coordinates one shows that $\omega_\sigma\wedge\omega=0$ and hence that $\omega\in\omega_\sigma\wedge\Omega^{p-1}(\A)=\Omega^{p-1}_0(\A)$ using \eqref{1}.
The claim follows by $\A$-normality of the latter module. 
\end{proof}

From Propositions~\ref{prp:proj} and \ref{prp:dual} we deduce the following dual version of Propositions~\ref{prp:proj}.
As in Definition~\ref{dfn:Omega-proj}, we define the sheaf $D_\bullet(\PP\A)$ on $\PP V$ by \eqref{def-Dp} in charts.

\begin{prp}\label{prp:proj-D}
As sheaves on $\PP V$, 
\[
\widetilde{D^0_p(\A)}=D_p(\PP\A).
\]
\end{prp}

\subsection{Self-duality of the relative log complex}

Let $dx=dx_1\wedge\cdots\wedge dx_\ell$. 
Then $\Omega^\ell(\A)=R\frac1f dx\cong R(n-\ell)$ as a graded $R$-module. 
On the other hand, contraction on $dx/f$ gives a graded map 
\begin{equation}\label{34}
D(\A)\to\Omega^{\ell-1}(\A)(n-\ell),\quad\delta\mapsto\iota_\delta(dx/f),
\end{equation}
and this can easily be seen to be an isomorphism. 
Let $\nu$ denote the image of the Euler differential:
\begin{equation}\label{eq:defnu}
\nu=\iota_\chi(dx/f)=\ell dx/df\in\Omega^{\ell-1}_\sigma(\A).
\end{equation}
By split exactness of $(\Omega^\bullet(\A),\omega_\sigma)$, or by \eqref{3} and \eqref{1}, 
\begin{equation}\label{9}
\omega_\sigma\colon\Omega^{\ell-1}_\sigma(\A)\to\Omega^\ell(\A)
\end{equation}
is an isomorphism with inverse $\iota_\chi$, and hence $\Omega^{\ell-1}_\sigma(\A)$ is a free rank one $R$-module generated by $\nu$.
By the same reason, 
\begin{equation}\label{20}
\Omega^{\ell-1}(\A)=\omega_\sigma\wedge\Omega^{\ell-2}_\sigma(\A)\oplus\Omega^{\ell-1}_\sigma(\A)\cong\Omega^{\ell-2}_\sigma(\A)\oplus\Omega^{\ell-1}_\sigma(\A).
\end{equation}
For $\delta\in D^\sigma(\A)$, $\iota_\delta(\omega_\sigma)=\delta(\alpha^\sigma)/\alpha^\sigma=0$ and hence
\[
\omega_\sigma\wedge\iota_\delta(\Omega^\ell(\A))=-\iota_\delta\left(\omega_\sigma\wedge\Omega^\ell(\A)\right)=0.
\]
Thus, $\iota_\delta(\Omega^\ell(\A))$ lies in the first summand of \eqref{20}.
By definition of the generator $\nu$ of $\Omega^{\ell-1}_\sigma(\A)$, $\iota_\chi(\Omega^\ell(\A))$ lies in the second summand of \eqref{20}.
Thus, composing \eqref{9} with $\iota_\delta$, the isomorphism \eqref{34}, induces an isomorphism
\begin{equation}\label{5}
D^\sigma(\A)\to\Omega^{\ell-2}_\sigma(\A)(\ell-n),\quad\delta\mapsto
\frac{1}{\omega_\sigma}\iota_\delta(\omega_\sigma\wedge\nu)=\iota_\delta(\nu).
\end{equation}
This proves, using Proposition~\ref{prp:dual} for the last part,

\begin{prp}\label{prp:self-dual}\
\begin{enumerate}
\item\label{sd1} $\Omega^\ell_\sigma(\A)=0$
\item\label{sd2} $\Omega^{\ell-1}_\sigma(\A)=R\nu\cong R(n-\ell)$
\item\label{sd3} $\Omega^{\ell-2}_\sigma(\A)\cong D^\sigma(\A)(n-\ell)\cong(\Omega_\sigma^1(\A)(\ell-n))^\vee$
\end{enumerate}
\end{prp}

Generalizing this result, consider non-degenerate pairing
\[
-\wedge-\colon\Omega_{V,f}^p\times\Omega_{V,f}^q\to
\Omega_{V,f}^\ell=R_f,\quad p+q=\ell,
\]
defined by the exterior product.
It induces a non-degenerate pairing
\begin{equation}\label{14}
-\wedge- \colon\Omega^p(\A)\times\Omega^q(\A)\to\Omega^\ell(\A)\cong R(n-\ell),\quad p+q=\ell,
\end{equation}
On $U$ this is easy to check in local coordinates using \eqref{13}, then it follows on $V$ by reflexivity.
It is immediate from the definition of the pairing that it turns $(\Omega^\bullet(\A),\omega)$ into a self-dual complex for any $\omega\in\Omega^1(\A)$.

\begin{prp}\label{prp:cplxdual}
The pairing \eqref{14} induces a non-degenerate pairing
\begin{equation}\label{15}
-\wedge-\colon\Omega^p_0(\A)\times\Omega^q_0(\A)\to\Omega_0^{\ell-1}(\A)\cong R(n-\ell),\quad p+q=\ell-1
\end{equation}
which makes the complex $(\Omega_0^\bullet(\A),\omega_\lambda)$ self-dual.
\end{prp}

\begin{proof}
Similar to \eqref{5}, \eqref{15} is induced by \eqref{14} by identifying $\Omega^i_0(\A)=\omega_\sigma\wedge\Omega^i(\A)\subset\Omega^{i+1}(\A)$ for $i=q,\ell$, and $\Omega^p_0(\A)=\Omega_\sigma^p(\A)$.
Well-definedness is then clear.
We have to show that for $\omega\in\Omega_0^q(\A)$ with $\Omega^p(\A)\wedge\omega_\sigma\wedge\omega=0$ we have $\omega_\sigma\wedge\omega=0$, or equivalently $\omega\in\omega_\sigma\wedge\Omega^{p-1}(\A)$.
But this hypothesis implies $\omega_\sigma\wedge\omega=0$ on $U$ and hence on $V$ by $\A$-normality of $\Omega^{q+1}(\A)$, so the claim follows.
\end{proof}

\section{Complexes}

\subsection{Parametric relative log complex}\label{ss:parametric}

The variety 
\begin{equation}\label{35}
\Sigma(\A)=\set{(x,\lambda)\in \big(V\setminus\bigcup_{H\in\A}H\big)\times\C^\A\mid \omega_\lambda(x)=0}\subset V\times\C^\A,
\end{equation}
was used in \cite{CDFV09} to study critical sets of master functions associated with the arrangement $\A$.
In \cite[Thm~2.9]{CDFV09}, it was shown that its closure 
\begin{equation}\label{38}
\overline{\Sigma}(\A)=\Spec(S/I),
\end{equation}
where the ideal
\[
I=I(\A):=\ideal{D(\A),\oma}\subset S:=R\otimes_\C C,\quad C:=\Sym((\C^\A)^\vee)=\C[\aa],
\]
was defined in \cite[\S2.5]{CDFV09} using the pairing \eqref{8}.
For the moment, we ignore the grading in $C$ and continue to refer only to the $\ZZ$-grading on $R$. 
Note that 
\begin{equation}\label{27}
H^\ell(\Omega_{S/C}^\bullet(\A),\oma)=(S/I(\A))dx/f\cong(S/I(\A))(\ell-n),
\end{equation}
where $\Omega_{S/C}^\bullet(\A)=\Omega^\bullet(\A)\otimes_\C C$.
In the spirit of this paper, let 
\begin{align}
\nonumber I_\sigma&=I_\sigma(\A):=\ideal{D^\sigma(\A),\oma}\subset S,\\
\label{47}I_0&=I_0(\A):=\ideal{D^0(\A),\oma}\subset S_0:=R\otimes_\C C_0,\quad C_0:=C/\ideal{|\aa|}.
\end{align}
Note that, by \eqref{28},
\begin{equation}\label{36}
S/I=S/(I_\sigma+\ideal{|\aa|})\cong S_0/I_0.
\end{equation}
To obtain the analogue of \eqref{27} for $I_0$ satisfying the parametric version of \eqref{3}, we need to work with parameters in $C_0$.
We use the notations
\[
\Omega_{S/C,\sigma}^\bullet(\A)=\Omega_\sigma^\bullet(\A)\otimes_\C C,\quad 
\Omega_{S_0/C_0,0}^\bullet(\A)=\Omega_0^\bullet(\A)\otimes_\C C_0.
\]
Since $\oma\wedge\iota_\delta(\nu)=\iota_\delta(\oma)\nu$, it follows from \eqref{5} that
\begin{gather}
\label{29}H^{\ell-1}(\Omega_{S/C,\sigma}^\bullet(\A),\oma)=
(S/I_\sigma)\nu\cong(S/I_\sigma)(n-\ell),\\
\label{30}H^{\ell-1}(\Omega_{S_0/C_0,0}^\bullet(\A),\oma)=
(S_0/I_0)\nu\cong(S_0/I_0)(n-\ell).
\end{gather}

Now assume that $\A$ is tame, so that $(\Omega_{S/C}^\bullet(\A),\oma)$ has cohomology concentrated in degree $\ell$ by \cite[Thm.~3.5]{CDFV09}.
For the next argument, we need a {\em second grading} for which the natural map $R\to\Omega^1_V$ has degree $-1$ and the coefficients $\aa$ have degree $1$.
Degree shifts with respect to this grading are denoted by square brackets $[-]$.
While $\oma$ is still of degree $0$, $\omega_\sigma$ has degree $-1$ for the second grading.  For the remainder of this paragraph, the differential in each complex is given by multiplication by $\oma$.
Consider now the exact sequence of second-graded complexes,
\[
\xymat@C=12pt{
0\ar[r] & \omega_\sigma\wedge\Omega_{S/C}^{\bullet-1}(\A)\ar[r] & \Omega_{S/C}^\bullet(\A)\ar[r]^-{\omega_\sigma} & \omega_\sigma\wedge\Omega_{S/C}^\bullet(\A)[-1]\ar[r] & 0,
}
\]
and that the second-degrees of $\omega_\sigma\wedge\Omega_{S/C}^\bullet(\A)$ have a lower bound, one sees that the cohomology of $\omega_\sigma\wedge\Omega_{S/C}^\bullet(\A)$ is concentrated in degree $\ell-1$.  From the exact sequences of $S$-modules
\begin{gather*}
\xymat@C-1em{
0\ar[r] & \omega_\sigma\wedge\Omega_{S/C}^{\bullet-1}(\A)\ar[r] & \Omega_{S/C}^\bullet(\A)\ar[r] & \Omega_{S/C,\sigma}^\bullet(\A)\ar[r] & 0,}\\
\xymat@C-1em{
0\ar[r] & \Omega_{S/C,\sigma}^\bullet(\A)\ar[r]^-{|\aa|} & \Omega_{S/C,\sigma}^\bullet(\A)\ar[r]\ar[r] & \Omega_{S_0/C_0,0}^\bullet(\A)\ar[r] & 0,}
\end{gather*} 
we see $\Omega_{S_0/C_0,0}^\bullet(\A)$ has cohomology at most in degrees $\ell-2$ and $\ell-1$.  
In fact there is no cohomology in degree $\ell-2$ since $|\aa|$ is a non-zero divisor on \eqref{29}.
To see this, choose $\sigma=e_1$, so that $I_\sigma$ has generators independent of $a_1$. 
By a linear change of coordinates in $C$, we may hence replace $|\aa|$ by the coordinate $a_1$, which is then clearly a non-zero divisor on $S/I_\sigma$.
Together with \eqref{30}, this proves the following:

\begin{prp}\label{prp:res0}
For tame $\A$, the complex $(\Omega_{S_0/C_0,0}^\bullet(\A),\oma)$ resolves $(S_0/I_0(\A))(n-\ell)$.
\end{prp}

In order to give this result a projective interpretation, let $\Gamma:=\Spec C_0$, consider the image $\PP_\Gamma\A$ of $\PP\A$ under the projection
\[
\pi_\Gamma\colon\PP_\Gamma V:=\Proj S_0\onto\PP V,
\]
and use Proposition~\ref{prp:proj} to identify
\[
\widetilde{\Omega_{S_0/C_0,0}^p(\A)}
=\Omega_{\PP_\Gamma V/\Gamma}^p(\PP_\Gamma\A)
=\pi_\Gamma^*\Omega^p(\PP\A).
\]
Consider also the preimage $\A_\Gamma$ of $\A$ under affine projection $V_\Gamma:=V\times\Gamma\to V$.
As in Proposition~\ref{22}, we can now identify $(\widetilde{\Omega_{S_0/C_0,0}^\bullet(\A)},\oma)$ and $(\Omega_{V_\Gamma/\Gamma}^\bullet(\A_\Gamma^{\{x_i=1\}}),\omega_{\aa_{\hat i}})$ where $\aa_{\hat i}$ is obtained by deleting the $i$th component from $\aa$. 
Setting $C_{\hat i}:=\C[a_j\mid j\ne i]\cong C_0$ and $\Gamma_{\hat i}:=\Spec C_{\hat i}$, the latter can then be identified with $(\Omega_{V_{\Gamma_{\hat i}}/\Gamma_{\hat i}}^\bullet(\A_\Gamma^{\{x_i=1\}}),\omega_{\aa_{\hat i}})$.

Similarly, $I_0$ can be related to 
\begin{equation}\label{48}
I(\PP\A):=\ideal{D(\PP\A),\oma}\subset\O_{\PP_\Gamma V}.
\end{equation}
Namely, by Proposition~\ref{prp:proj-D}, we have that 
\begin{equation}\label{37}
\widetilde{I_0(\A)}=I(\PP\A).
\end{equation}
A dual version of \eqref{21} serves to identify $I_0(\A)_{x_i}=I_{e_i}(\A)_{x_i}$ and $I(\A^{\{x_i=1\}})$, where the latter is the restriction of the ideal sheaf $I(\PP\A)$ to the chart $\{x_i=1\}$.
Then we can also relate $\overline{\Sigma}(\A)$ to 
\[
\overline\Sigma(\PP\A):=\Proj(S_0/I(\PP\A)).
\]
Recall from \eqref{35} that $\overline{\Sigma}(\A)$ is a subvariety of $V\times\C^\A$; 
let $\PP\overline{\Sigma}(\A)$ denote its projectivization in the first factor.
By \eqref{38}, \eqref{36}, and \eqref{37}, $\overline\Sigma(\PP\A)=\PP\overline{\Sigma}(\A)$ via the natural inclusion $\PP_\Gamma V\into\PP V\times \C^\A$.

Finally, the projective version of Proposition~\ref{prp:res0} reads

\begin{prp}\label{prp:res0-proj}
For tame $\A$, $(\pi_\Gamma^*\Omega^\bullet(\PP\A),\oma)$ resolves $\O_{\overline{\Sigma}(\PP\A)}(n-\ell)$.
\end{prp}

\subsection{Log complex in tame and generic cases}

Specializing \eqref{47} and \eqref{48} to $\aa=\lambda\ne0$, $|\lambda|=0$, we set (as in \cite[Prop.~2.7]{OT95}):
\begin{align*}
I_{0,\lambda}=I_{0,\lambda}(\A)&:=\ideal{D^0(\A),\omega_\lambda}\subset R,\\
I_\lambda(\PP\A)&:=\ideal{D(\PP\A),\omega_\lambda}\subset\O_{\PP V},
\end{align*}
such that \eqref{37} gives
\begin{equation}\label{49}
\widetilde{I_{0,\lambda}(\A)}=I_\lambda(\PP\A).
\end{equation}
We denote by $\overline\Sigma_\lambda(\A)$ and $\overline\Sigma_\lambda(\PP\A)$ the specializations of $\overline\Sigma(\A)$ and $\overline\Sigma(\PP\A)$ to $\aa=\lambda$ respectively.

Proposition~\ref{prp:res0} has the following application, which improves \cite[Prop.~3.9]{CDFV09}:

\begin{prp}\label{prp:critset}
If $\A$ is a tame arrangement, then  $H^p(\Omega^\bullet_0(\A),\omega_\lambda)\neq0$ implies that
the codimension of $\overline\Sigma_\lambda(\A)$ is at most $p$, provided that either $\A$ has rank at most $4$, $\A$ is free, or $p\leq 2$.
\end{prp}

\begin{proof}
For $\lambda\in\C^\A$, let $R_\lambda=S/\ideal{\set{a_i-\lambda_i}}$.
Arguing as in \cite[Prop.~3.9]{CDFV09}, Proposition~\ref{prp:res0} implies the hyper-$\Ext$ spectral sequence
\begin{equation}\label{eq:hypExt}
E^{p,q}_1=\Ext^q_{S_0}(\Omega^{\ell-1-p}_{S_0/C_0,0}(\A),R_\lambda)
\cong\Ext^q_R(\Omega^{\ell-1-p}_0(\A),R)
\end{equation}
converges to $\Ext^{p+q}_{S_0}(S_0/I_0(\A),R_\lambda)$, and $E^{p,0}_1\cong\Omega^p_0(\A)$ by Proposition~\ref{prp:cplxdual}.
Since $E^{p0}_2$ receives no nonzero differentials for $p\leq 2$, we see that, if $H^p(\Omega^\bullet_{S_0/C_0,0}(\A),\omega_\lambda)$ is nonzero, then so is
$\Ext^p_{S_0}(S_0/I_0(\A),R_\lambda)$ for $p\leq2$, as well as for $p=3$ when $\ell-1=3$. 
This implies that the codimension of $\overline\Sigma_\lambda(\A)$ is at most $p$.
\end{proof}

The projective version of this result based on Proposition~\ref{prp:res0-proj} reads:

\begin{prp}
If $\A$ is a tame arrangement, then  $H^p(\pi_\Gamma^*\Omega^\bullet(\PP\A),\omega_\lambda)\neq0$ implies that the codimension of $\overline{\Sigma}_\lambda(\PP\A)$ is at most $p$, provided that either $\A$ has rank at most $4$, $\A$ is free, or $p\leq 2$.
\end{prp}

In \cite{OT95}, Orlik and Terao show that there is a Zariski-open subset $Y\subseteq\C^\A$ with the property that, for $\lambda\in Y$, the function \eqref{eq:master} has non-degenerate, isolated critical points.  We will call such $\lambda$ {\em generic}.  
Part of their argument shows that, for a noncentral, affine arrangement $\A$, the cohomology of the complex $(\Omega^\bullet(\A),\omega_\lambda)$ is concentrated in top dimension, for $\lambda\in Y$: see \cite[Prop~4.6]{OT95}.
We can prove a projective version of this result, which is a slightly stronger statement, but requires a tameness hypothesis.

\begin{dfn}
We call $\A$ {\em almost tame} if $\pd\Omega^p(\A)\le p+1$ for $p=1,\dots,\ell$.
\end{dfn}

As in Definition~\ref{df:tame}, the condition is vacuous except for $1\le p\le\ell-4$. 
For the sake of stating some results as generally as possible, we have now introduced two new homological boundedness conditions on an arrangement. 
For the reader's convenience, we summarize their relative strength as follows:
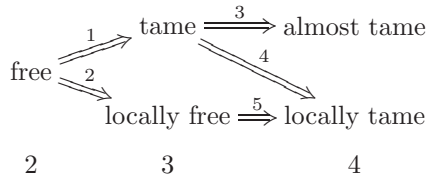
\[
\xymat@C-1em@R-2em{
&\text{tame}\ar@{=>}[r]^-{3}\ar@{=>}[ddr]^-{4} & \text{almost tame}\\
\text{free}\ar@{=>}[ur]^-{1}\ar@{=>}[dr]^-{2}\\
&\text{locally free}\ar@{=>}[r]^-{5} & \text{locally tame}\\ 
2 & 3 & 4
}
\]
The properties in the column labeled $\ell$ hold for all arrangements of rank $\leq\ell$. 
\cite[Example~5.3]{CDFV09} gives an arrangement of rank $4$ which is locally free but not tame, showing implications $2$, $3$ and $4$ cannot be reversed in general. 
The rank-$4$ arrangements of Example~\ref{ex:ziegler} are tame but not locally free, so implications $1$ and $5$ also cannot be reversed.

\begin{prp}\label{prp:genExact}
If $\A$ is almost tame and $\lambda$ is generic (with $|\lambda|=0$), then 
\begin{equation}\label{eq:genExact}
H^p(\Omega^\bullet_0(\A),\omega_\lambda)\cong
\begin{cases}
0&\text{for }p\ne\ell-1,\\
(R/I_{0,\lambda}(\A))(n-\ell)&\text{for }p=\ell-1.
\end{cases}
\end{equation}
\end{prp}

\begin{proof}
Denote by $\hat\Omega^\bullet_0(\A)$ the complex obtained from $(\Omega^\bullet_0(\A),\omega_\lambda)$ by replacing $\Omega^{\ell-1}_0(\A)$ by $\omega_\lambda\wedge\Omega^{\ell-2}_0(\A)$.
Consider the two hyperhomology spectral sequences of the complex $\hat\Omega^\bullet_0(\A)$:
\begin{equation}\label{2}
^IE_{p,q}^2=\Ext_R^{-q}(H^p(\hat\Omega^\bullet_0(\A),R),\quad 
^{II}E_{p,q}^1=\Ext_R^{-q}(\hat\Omega^p_0(\A),R).
\end{equation}
By Proposition~\ref{22} and \cite{OT95}, it follows that $\hat\Omega^\bullet_0(\A)$ is exact away from the origin.
So $^IE_{p,-q}^2$ is non-zero only for $q=\ell$ and $p<\ell-1$, contributing to degrees $p-q<-1$ in the abutment.
But by hypothesis and \eqref{1}, $^{II}E_{p,-q}^1$ is non-zero only if $q\le p+1$, contributing to degrees complementary to those of $^IE_{p,q}^2$.
This shows that both sequences converge to zero.
But the first sequence degenerates on the $E^2$-page, and hence $^IE_{p,q}^2=0$.
This proves the first claim; the second follows by specializing \eqref{30} to $\aa=\lambda$.
\end{proof}

The projective version of this result does not require tameness, as we can apply \cite{OT95} in charts.

\begin{prp}\label{prp:genExact-proj}
If  $\lambda$ is generic (with $|\lambda|=0$), then 
\begin{equation}\label{eq:genExact-proj}
H^p(\Omega^\bullet(\PP\A),\omega_\lambda)\cong
\begin{cases}
0,&p\ne\ell-1,\\
(\O_{\PP V}/I_\lambda(\PP\A))(n-\ell),&p=\ell-1.
\end{cases}
\end{equation}
\end{prp}

In special cases, a more detailed understanding is possible.  For example, from \cite[Thm.~5]{CHKS06} we see that if $\A$ is a generic arrangement, the critical set of the master function \eqref{eq:master} in $\PP V$ is zero-dimensional, for all nonzero $\lambda$ with $\abs{\lambda}=0$. 
From this we note the following:

\begin{prp}
If $\A$ is a generic arrangement, then \eqref{eq:genExact} holds for all $\lambda\ne0$ (with $\abs{\lambda}=0$).
\end{prp}

\begin{proof}
In \cite{DSSWW10}, we show that if $\A$ is generic, then $\Ext^q_R(\Omega^p_0(\A),R)=0$ except for $q=0$ and $q=p$.  Then the $E_2$-page of the spectral sequence \eqref{eq:hypExt} is zero except for $p+q=\ell-1$ and $q>0$, while
$E_2^{p,0}=H^p(\Omega^p_0(\A),\omega_\lambda)$ by Proposition~\ref{prp:cplxdual}.  

However, by \cite[Thm.~5]{CHKS06}, the codimension of $\overline{\Sigma}_\lambda(\PP\A)$, and hence of $\overline{\Sigma}_\lambda(\A)$, is $\ell-1$, so $E_\infty^{p,q}=0$ for 
$p+q\leq \ell-1$.  It follows that $E_2^{p,0}=0$ for $0\leq p\leq \ell-2$.
\end{proof}

Passing to coherent sheaves, we obtain the projective analogue of Proposition~\ref{prp:genExact-proj} as well.

\begin{cor}
If $\A$ is a generic arrangement, then \eqref{eq:genExact-proj} holds for all $\lambda\ne0$ (with $\abs{\lambda}=0$).
\end{cor}

\section{Chern classes}\label{sec:chern}

In this section, we prove an analogue of the Borel--Serre formula for sheaves on $\PP^d$ with projective resolution of length one (see Theorem~\ref{thm:BorSer}).
Then we apply this formula and the theory developed in the preceding sections to prove a generalized Musta{\c{t}}{\u{a}}--Schenck formula for tame arrangements with zero-dimensional non-free locus (see Theorem~\ref{th:Chern}).

\subsection{Polynomial identities}

We begin with some technical preparations for the following sections.
First, we work in the ring $\Q[[u,t]]$.
Consider the power series
\begin{equation}\label{eq:defF}
F_\gamma(t,u)=\prod_{i=1}^n(1+ue^{\gamma_i t})\in\Q[[u,t]],
\end{equation}
with parameters in the ring of symmetric functions in a set of variables $\gamma$, for which we refer to \cite{MacBook}.
Let $\alpha=\set{\alpha_1,\ldots,\alpha_n}$ and $\beta=\set{\beta_1,\ldots,\beta_{n-r}}$ denote two sets of variables, and let
\begin{equation}\label{eq:defC}
C(t,u)=F_\alpha(t,u)/F_\beta(t,u)\in\Q[[u,t]].
\end{equation}
Let 
\begin{align*}
E_\gamma(t)&=\prod_{i=1}^n (1+\gamma_i t)\in\Q[t],\\
H_\gamma(t)&=\prod_{i=1}^n (1-\gamma_i t)^{-1}\in\Q[[t]]
\end{align*}
denote the generating function for the elementary symmetric and complete symmetric functions respectively, where the variables are $\gamma=\alpha$ or $\gamma=\beta$.

Denote by $\SS\subset\Q[[u,t]]$ the subset of power series in $u$ and $t$, for which the coefficient of $t^k$ is a polynomial in $u$ of degree at most $k$. 
It is easy to see that $\SS$ is closed under taking products and multiplicative inverses (whenever defined).

\begin{lem}\label{lem:identity}
$\F_\gamma((1+u)t,u)/(1+u)^n\in\SS$ for any variables $\gamma=\set{\gamma_1,\ldots,\gamma_n}$, and 
\[
\left.\frac{F_\gamma((1+u)t,u)}{(1+u)^n}\right|_{u=-1}=E_\gamma(-t).
\]
\end{lem}

\begin{proof}
We expand $F_\gamma((1+u)t,u)$ as a power series in $t$. 
Note that, for any $i$,
\begin{align*}
1+ue^{\gamma_i(1+u)t}
&=1+u+\gamma_iu(1+u)t+u(1+u)^2t^2P_i((1+u)t)\\
&=(1+u)(1+\gamma_iut+u(1+u)t^2P_i((1+u)t)),
\end{align*}
for some power series $P_i(t)$, and hence
\[
F_\gamma((1+u)t,u)/(1+u)^n=\prod_{i=1}^n(1+\gamma_i u t +u(u+1)t^2P_i((u+1)t)),
\]
from which we see the coefficient of $t^k$ is a polynomial of degree $k$ in $u$, so we may evaluate to obtain
\[
\left.F_\gamma((1+u)t,u)/(1+u)^n\right|_{u=-1}=E_\gamma(-t),
\]
as required.
\end{proof}

\begin{prp}\label{prp:identity}
One can write
\[
C(t,u)=\sum_{k\geq0}(1+u)^{r-k}a_k(u)t^k,
\]
where $a_k(u)$ are polynomials in $u$ of degree at most $k$, such that
\begin{equation}\label{eq:rephrase}
\sum_{k\geq0}a_k(-1)t^k=E_\alpha(-t)H_\beta(t).
\end{equation}
\end{prp}

\begin{proof}
By Lemma~\ref{lem:identity}, $F_\alpha((1+u)t,u)/(1+u)^n$ and $F_\beta((1+u)t,u)/(1+u)^{n-r}$ are both in $\SS$, so $C((1+u)t,u)/(1+u)^r$ is too. 
Moreover,
\begin{eqnarray*}
\left.\frac{C((1+u)t,u)}{(1+u)^{r}}\right|_{u=-1} &=& E_\alpha(-t)/E_\beta(-t)\\
&=& E_\alpha(-t)H_\beta(t).
\end{eqnarray*}
The claim follows.
\end{proof}

Now we switch to the ring $A[[u]]$, where
\[
A:=\Q[t]/\ideal{t^{d+1}},
\]
and consider the image $\overline{C}(t,u)$ in $A[[u]]$ of $C(t,u)$ from \eqref{eq:defC}.
If $r\geq d$, then,  by Proposition~\ref{prp:identity}, $\overline{C}(t,u)$ becomes a polynomial in $u$ of degree $r$, which motivates the following. 

\begin{dfn}\label{dfn:Lpoly}
For $0\leq p\leq r\ge d$, we call $L^p_{\alpha,\beta}(t)\in\Q[t]/\ideal{t^{d+1}}$ defined by 
\begin{equation}\label{eq:identity2}
\sum_{p=0}^rL^p_{\alpha,\beta}(t)u^p:=\overline{C}(t,u).
\end{equation}
the \emph{$p$th Lebelt polynomial}.
\end{dfn}

For any variables $\gamma$, let $\abs{\gamma}=\sum_{i=1}^n\gamma_i$.

\begin{lem}\label{lem:chtop}
For $r=d$, 
\begin{align*}
L^r_{\alpha,\beta}&=e^{(\abs{\alpha}-\abs{\beta})t},\text{~and}\\
L^{r-1}_{\alpha,\beta}&=e^{(\abs{\alpha}-\abs{\beta})t}\left(
\sum_{i=1}^ne^{-\alpha_i t}-
\sum_{i=1}^{n-r} e^{-\beta_i t}\right).
\end{align*}
\end{lem}

\begin{proof}
As $\overline{C}(t,u)$ is a polynomial in $u$ of degree $r$, we may change variables to get
\[
L^r_{\alpha,\beta}\overline{C}(t,u)
=\left.u^r\overline{C}(t,u^{-1})\right|_{u=0}
=\left.\frac
{\prod_{i=1}^{n}(u+e^{\alpha_i t})}
{\prod_{i=1}^{n-r}(u+e^{\beta_i t})}\right|_{u=0}
=\frac{\prod_{i=1}^ne^{\alpha_i t}}{\prod_{i=1}^{n-r}e^{\beta_i t}},
\]
as required.
For the second claim, 
\[
L^{r-1}_{\alpha,\beta}
=\frac{d}{du}\left.\frac
{\prod_{i=1}^n(u+e^{\alpha_i t})}
{\prod_{i=1}^{n-r}(u+e^{\beta_i t})}\right|_{u=0}\\
=e^{(\abs{\alpha}-\abs{\beta})t}\left(
\sum_{i=1}^ne^{-\alpha_i t}-
\sum_{i=1}^{n-r} e^{-\beta_i t}\right).
\]

\end{proof}

\subsection{Lebelt resolutions}\label{sec:Lebelt}

We shall now consider a coherent sheaf of rank $r$ on $\PP^d$ having a projective resolution of length one,
\begin{equation}\label{eq:resM}
\xymat@C-0.5em{
0\ar[r] & \F_1\ar[r] & \F_0\ar[r] & \M\ar[r] & 0.
}
\end{equation}
For such sheaves we shall prove a Borel--Serre formula in Theorem~\ref{thm:BorSer}.
We will compute in the extended rational Chow ring
\[
A(\PP^d)_\Q[[u]]\cong A[[u]].
\]
replacing $\alpha$ and $\beta$ by the Chern roots of $\F_0$ and $\F_1$,

\begin{prp}\label{prp:wedges}
Assume that $\M$ is locally a $(k-1)$st syzygy. 
Then 
\[
\ch_t\bigwedge^p\M=L^p_{\alpha,\beta}
\]
for $0\leq p\leq k$. 
In particular, this holds vacuously for $p=1$.
\end{prp}

\begin{proof}
By the hypothesis on $\M$, we may use the Lebelt resolution~\cite{Leb77} 
\[
\xymat@C-15pt{
0\ar[r] & S^p\F_1\ar[r] & S^{p-1}\F_1\otimes\bigwedge^1\F_0\ar[r] & \cdots\ar[r] & S^1\F_1\otimes\bigwedge^{p-1}\F_0\ar[r] & \bigwedge^p\F_0\ar[r] & \bigwedge^p\M\ar[r] & 0
}
\]
to resolve $\bigwedge^p\M$, for all $p\leq k$.
Here we are using the fact that local resolutions glue, by the uniqueness of the Lebelt's differential proved in \cite[(2.2)~Satz]{Leb77}.
Then 
\begin{align*}
\ch_t(\bigwedge^s\M)u^s
&=\sum_{p+q=s}(-1)^q\ch_t(S^q\F_1)\ch_t(\bigwedge^p\F_0)\\
&=\sum_{p+q=s}(-1)^q
\sum_{i_1<\cdots<i_p}ue^{\alpha_{i_1}t}\cdots ue^{\alpha_{i_p}t}
\sum_{j_1\leq \cdots\leq j_q}ue^{\beta_{j_1}t}\cdots ue^{\beta_{j_q}t},
\end{align*}
using the splitting principle for symmetric and exterior powers of vector bundles.
The claim now follows by expanding the expression \eqref{eq:defC} as a power series in $u$, and noting that this is the $u^p$-term, which is $L^p_{\alpha,\beta}$ by Definition~\ref{dfn:Lpoly}.
\end{proof}

The equation \eqref{eq:rephrase} gives the Chern polynomial of $\M$:

\begin{lem}\label{lem:ct}
\[
c_t(\M)=\sum_{k=0}^d a_k(-1)(-t)^k.
\]
\end{lem}

\begin{proof}
By \eqref{eq:resM} and the multiplicativity of Chern polynomials, 
\[
c_t(\M)=c_t(\F_0)/c_t(\F_1)=E_\alpha(t)/E_\beta(t)=E_\alpha(t)H_\beta(-t)
\]
and the result follows by \eqref{eq:rephrase}.
\end{proof}

For $r=d$, we can prove an analogue of the Borel--Serre formula for vector bundles, by expressing the top Chern class of $\M$ in terms of Lebelt polynomials.

\begin{thm}\label{thm:BorSer}
For $r=d$,
\[
c_{d}(\M)=(-1)^{d}\sum_{p=0}^{d}(-1)^p L^p_{\alpha,\beta}.
\]
\end{thm}

\begin{proof}
Since $\overline{C}(t,u)$ is a polynomial in $u$, we may set $u=-1$ in
\eqref{eq:identity2} to get
\[
\sum_{p=0}^d(-1)^p L^p_{\alpha,\beta}
=\overline{C}(t,-1)
=a_{d}(-1).
\]
Now by Lemma~\ref{lem:ct}, our
Euler characteristic is obtained as the coefficient of $t^d$.
\end{proof}

We can mimic another formula for vector bundles:
If $\M$ was locally free of rank $r=d$, the we would have an isomorphism $\bigwedge^{r-1}\M\cong\bigwedge^r\M\oplus\M^\vee$, and hence $\ch_t(\bigwedge^r\M)\cdot\ch_{-t}\M=\ch_t(\bigwedge^{r-1}\M)$.
Lemma~\ref{lem:chtop} proves the following analogue of this formula:

\begin{prp}\label{prp:chtop}
For $r=d$, 
\begin{align*}
e^{c_1(\M)t}&=L^r_{\alpha,\beta},\\
e^{c_1(\M)t}\ch_{-t}(\M)&=L^{r-1}_{\alpha,\beta}.
\end{align*}
\end{prp}

\subsection{Musta{\c{t}}{\u{a}}--Schenck formulas}

We recall a result of Musta{\c{t}}{\u{a}} and Schenck \cite[Thm.~4.1]{MuSc01}, which we formulate projectively here.
Due to a different grading convention, $\Omega^1(\A)$ in loc.\ cit.\ translates to 
\[
\Omega^1(\A)(1)\cong\Omega^1_0(\A)(1)\oplus S(1)
\]
and its associated total Chern polynomial in the Chow ring $A(\PP^{\ell-1})\cong\ZZ[t]/\ideal{t^\ell}$ to
\[
c_t(\Omega^1(\A)(1))=(1+t)c_t(\Omega^1(\PP\A)(1))
\]
using Proposition~\ref{prp:proj} and multiplicativity of $c_t$.
Denote by $\pi(\A,t)$ and $\pi(\PP\A,t)$ the Poincar\'e polynomials of the complements $V\setminus\A$ and $\PP V\setminus\PP\A$ respectively.
By \cite[Prop.~2.51]{OTbook} (see also \cite[Prop.~5.1, Thm.~5.90]{OTbook}) they are related by 
\[
\pi(\A,t)=(1+t)\pi(\PP\A,t).
\]
Now the projective version of \cite[Thm.~4.1]{MuSc01} reads

\begin{thm}[\cite{MuSc01}]\label{th:locfreeChern}
Let $\PP\A$ be a locally free arrangement in $\PP V$. 
Then
\begin{equation}\label{eq:locfreeChern}
c_t(\Omega^1(\PP\A)(1))=\pi(\PP\A,t)
\end{equation}
in the Chow ring $A(\PP V)\cong\ZZ[t]/\ideal{t^\ell}$ where $\pi(\PP\A,t)$ is the Poincar\'e polynomial of the complement $\PP V\setminus\PP\A$, and $c_t$ denotes the total Chern polynomial.
\end{thm}

In the following, we generalize this result for locally tame arrangements in $\PP^{\ell-1}$ with zero-dimensional non-free locus.
Using that $\A\mapsto\Omega^1(\A)$ is a local functor, we have
\begin{equation}\label{eq:projExt}
\SExt^p_{\O_{\PP V}}(\Omega^1(\PP\A),\O_{\PP^d})=
\begin{cases}
\bigoplus\limits_{X\in L_{\ell-1}(\A)}\Ext^1_{R_X}(\Omega^1(\A_X),R_X),&\text{ if }p=1,\\
0,&\text{ if }p\ge2.
\end{cases}
\end{equation}
where we consider $\A_X$ as an arrangement in an affine chart $\AA^{\ell-1}$ of $\PP V$ with origin $X$ and coordinate ring $R_X$.

\begin{dfn}\label{def:N}
For a central arrangement $\A$ in $V$ with zero-dimensional non-free locus, let 
\[
N(\A)=\length(\Ext^1_R(\Omega^1(\A),R),
\]
a non-negative integer. 
For an arrangement $\PP\A$ in $\PP V$ with zero-dimensional non-free locus, let 
\begin{equation}\label{eq:Nloc}
N(\PP\A)=h^0(\SExt^1_{\O_{\PP V}}(\Omega^1(\PP\A),\O_{\PP V})=\sum_{X\in L_{\ell-1}(\A)}N(\A_X),
\end{equation}
using \eqref{eq:projExt}.
\end{dfn}

Note that $N(\A)=0$ if and only if $\A$ is free, and that $N(\PP\A)=0$ if and only if $\PP\A$ is locally free.

Recall that a rank-$\ell$ arrangement $\A$ is called {\em $k$-generic} if $\#\A_X=k$ for all $X\in L_k(\A)$, and $\A$ is called generic if it is {\em $(\ell-1)$-generic}.

By \cite{DSSWW10}, $N(\A)={n-1\choose\ell}$ for a generic, rank-$\ell$ arrangement $\A$ of $n$ hyperplanes.
This gives the following explicit combinatorial formula:

\begin{thm}
If $\PP\A$ is a $\ell-2$-generic arrangement in $\PP^{\ell-1}$, then
\[
N(\PP\A)=\sum_{X\in L_\ell(\A)}{\#\A_X-1\choose\ell}.
\]
\end{thm}

If $\A$ is tame, one can compute $N(\A)$ by comparing the Hilbert series of logarithmic forms and differentials:

\begin{prp}\label{prp:N}
For a central arrangement $\A$ with a zero-dimensional non-free locus,
\begin{equation}\label{eq:N}
N(\A)=\lim_{t\to1} h(D_1(\A),t)-h(\Omega^1(\A),t^{-1})/(-t)^\ell,
\end{equation}
where $h(M,t)$ denotes the Hilbert series of a graded module $M$.
\end{prp}

\begin{proof}
By the tame hypothesis, $\pd\Omega^1(\A)\leq1$, so we have a graded free resolution of the form
\[
\xymat@C-0.5em{
0 & \Omega^1(\A)\ar[l] & F_0\ar[l] & F_1\ar[l] & 0\ar[l].
}
\]
Dualizing gives
\[
\xymat@C-0.5em{
0\ar[r] & D_1(\A)\ar[r] & F_0^\vee\ar[r] & F_1^\vee\ar[r] & \Ext^1_R(\Omega^1(\A),R)\ar[r] & 0.
}
\]
Then $h(F^\vee,t)=h(F,t^{-1})/(-t)^\ell$ if $F$ is a free module over a graded ring of dimension $\ell$, so
\begin{align*}
h(\Ext^1_R(\Omega^1(\A),R),t)
&=h(F_1^\vee,t)-h(F_0^\vee,t)+h(D_1(\A),t)\\
&=-h(\Omega^1(\A),t^{-1})/(-t)^\ell+h(D_1(\A),t).
\end{align*}
The right-hand side, then, reduces to a polynomial in $t$, so we may compute $N(\A)$ by evaluating at $t=1$, giving \eqref{eq:N}.
\end{proof}

We can now measure the extent to which the formula of Theorem~\ref{th:locfreeChern} fails for locally tame arrangements with zero-dimensional non-free locus.

\begin{thm}\label{th:Chern}
If $\PP\A$ is a locally tame arrangement in $\PP V$ with zero-dimensional non-free locus, then
\begin{equation}\label{eq:Chern}
c_t(\Omega^1(\PP\A)(1))=\pi(\PP\A,t)+N(\PP\A)t^{\ell-1}.
\end{equation}
\end{thm}

For $\ell=4$, both hypotheses of Theorem~\ref{th:Chern} are trivially fulfilled by reflexivity of $\Omega^1(\PP\A)$ and \cite[Prop.~1.3, Cor.~1.4]{Har80}.

\begin{cor}
If $\PP\A$ is an arrangement in $\PP^3$, then
\[
c_t(\Omega^1(\PP\A)(1))=\pi(\PP\A,t)+N(\PP\A)t^3.
\]
\end{cor}

\begin{cor}
If $\PP\A$ is an arrangement in $\PP^3$, or a locally tame arrangement in $\PP V$ with zero-dimensional non-free locus, the Musta{\c{t}}{\u{a}}--Schenck formula \eqref{eq:locfreeChern} holds if and only if $\A$ is locally free.
\end{cor}

For an arrangement $\PP\A$ in $\PP V$, we denote
\begin{align}
\label{40}c_t(\Omega^1(\PP\A)(1))&=1+c_1t+c_2t^2+\cdots+c_{\ell-1}t^{\ell-1},
\quad\text{and}\\
\label{41}\pi(\PP\A,t)&=1+b_1t+b_2t^2+\cdots+b_{\ell-1}t^{\ell-1}.
\end{align}
Since $c_1c_2\equiv c_3\mod2$ for Chern classes of coherent sheaves on $\PP^3$ (see, e.g., \cite[Cor.~2.4]{Har80}), we observe:

\begin{cor}
If $\PP\A$ is an arrangement in $\PP^3$, then 
\[
N(\PP\A)\equiv b_1b_2+b_3\mod2,
\]
where the $b_1$, $b_2$, and $b_3$ are the coefficients of the Poincar\'e polynomial \eqref{41}.
In particular, if $\PP\A$ is locally free, then $b_1b_2\equiv b_3\mod2$.
\end{cor}

The proof will follow some preliminary observations for general $\ell$:
Let $H_0$ be a hyperplane in a rank-$\ell$ arrangement $\A$, and let $\A'=\A\setminus\set{H_0}$. 
Recall that $\A^{H_0}=\set{H\cap H_0\mid H\in\A'}$, an arrangement in $H_0$. 
The inclusion of $H_0\subset V$ gives a map $i\colon\PP H_0\to\PP V$. 
By \cite[Cor.~4.5]{Zi89}, restriction gives
an exact sequence:
\begin{equation}\label{eq:omegas}
\xymat{
0\ar[r] & \Omega^p(\A)(-1)\ar[r]^-{\alpha_{H_0}} & \Omega^p(\A')\ar[r]^-{i^*} & \Omega^p(\A^{H_0}).
}
\end{equation}
The hyperplane $H_0\in\A$ is called {\em generic} if $\codim_{H_0} H_0\cap X=\codim X$ for all $X\in L_{<\ell}(\A')$. 
Ziegler notes in \cite[Ex.~8.7.(iii)]{Zi89} that $i^*$ need not be surjective, even when $H_0$ is generic. 

However, this subtlety happens at the origin and therefore disappears if one passes to projective space.
Note that taking the kernel of $\iota_\chi$ in \eqref{eq:omegas} and then sheafifying yields, using Proposition~\ref{prp:proj}, the sequence \eqref{70} below, and its exactness except at the right.

\begin{prp}\label{prp:slice-res}
Let $\PP\A$ be a projective arrangement, and $H_0\in\A$ generic.
Then there are exact sequences of sheaves on $\PP V$ and $\PP H_0$
respectively,
\begin{gather}
\label{70}\xymat{
0\ar[r] & \Omega^p(\PP\A)(-1)\ar[r]^-{\alpha_{H_0}} &
\Omega^p(\PP\A')\ar[r]^-{i^*} & i_*\Omega^p(\PP\A^{H_0})\ar[r] & 0,
}\\
\label{71}\xymat{
0\ar[r] & \Omega^{p-1}(\PP\A^{H_0})(-1)\ar[r]^-{d\alpha_{H_0}} &
i^*\Omega^p(\PP\A')\ar[r]^-{i^*} & \Omega^p(\PP\A^{H_0})\ar[r] & 0,
}
\end{gather}
for $0\leq p\leq\ell-1$.
\end{prp}

\begin{proof}
Genericity of $H_0$ implies that, for any closed point $x\in\PP
H_0\cap\PP\A'$, $X=\bigcap_{x\in\PP H\in\PP\A'}\PP H$ is at least
$1$-dimensional, and $\PP\A$ is, locally at $x$, a product of $\PP\A'_X$
with the hyperplane $\PP H_0$.
Thus,
\begin{gather*}
\Omega^p(\PP\A')_x=\pi^*\Omega^p(\PP\A^{H_0})_x\oplus
d\alpha_{H_0}\wedge\pi^*\Omega^{p-1}(\PP\A^{H_0})_x,\\
\Omega^p(\PP\A)_x=\pi^*\Omega^p(\PP\A^{H_0})_x\oplus
\frac{d\alpha_{H_0}}{\alpha_{H_0}}\wedge\pi^*\Omega^{p-1}(\PP\A^{H_0})_x,
\end{gather*}
where $\pi\colon\PP V\setminus\{\alpha_{H_0}=\infty\}\to\PP H_0$ is a
projection.
The claim follows.
\end{proof}

\begin{prp}\label{prp:slice-Chern}
Let $\PP\A$ be a projective arrangement whose non-free locus has codimension $>k$.
Then the polynomials \eqref{40} and \eqref{41} agree through degree $k$, that is, $c_i=b_i$ for $i=1,\dots,k$.
In particular, this holds true for $k=2$ without any hypothesis.
\end{prp}

\begin{proof}
We argue by induction on $\ell$, and we may assume that $k\ge2$ by reflexivity of $\Omega^1(\PP\A)$ and  \cite[Cor.~1.4]{Har77}.
If $\ell\leq k+1$ then $\A$ is locally free, and the result follows by Theorem~\ref{th:locfreeChern}. 
Otherwise, let $H_0\subset V$ be a linear hyperplane that meets $\A$ generically. 
Then $L_{\leq \ell-2}(\A)=L_{\leq \ell-2}(\A^{H_0})$, hence
\[
\pi(\PP\A^{H_0},t)=1+b_1t+\cdots+b_{l-2}t^{l-2}.
\]
We obtain a short exact sequence from \eqref{71} of Proposition~\ref{prp:slice-res}:
\[
\xymat@C-0.5em{
0\ar[r] & \O_{H_0}\ar[r] & i^*\Omega^1(\PP\A)(1)\ar[r] & \Omega^1(\PP\A^{H_0})(1)\ar[r] & 0,
}
\] 
from which we see that $c_t(i^*\Omega^1(\PP\A)(1))=c_t(\Omega^1(\PP\A^{H_0})(1))$.

Now the codimension of the non-free locus in $\A^{H_0}$ is at least
as large as that in $\A$, so by induction,
\[
c_t(\Omega^1(\PP\A^{H_0})(1))\equiv 1+b_1t+\cdots+b_kt^k\mod t^{k+1}.
\]

Now let $i\colon\PP {H_0}\into\PP V$ denote the inclusion map; the induced map on Chow rings agrees with the canonical map 
\[
i^*\colon\ZZ[t]/t^{\ell}\to\ZZ[t]/t^{\ell-1}.
\]
Using the functoriality of Chern classes (see, e.g., \cite[\S A.3]{Har77}) for the third equality, we deduce
\begin{align*}
1+b_1t+\cdots+b_kt^k
&\equiv c_t(\Omega^1(\PP\A^{H_0})(1))\\
&=c_t(i^*\Omega^1(\PP\A)(1))\\
&=i^*c_t(\Omega^1(\PP\A)(1))\\
&\equiv 1+c_1t+\cdots+c_kt^k\mod t^{k+1}.
\end{align*}
\end{proof}

\begin{lem}\label{lem:ch-pt}
The Chern character of a reduced point in $\PP^d$ is $t^d$.
\end{lem}
\begin{proof}
The point has a Koszul resolution $\bigwedge^\bullet\E$ where $\E=\O_{\PP^d}(-1)^d$.
By \cite[Ex.~3.2.5]{FultonBook}, its Chern character equals $c_d(\E^\vee)t^d/\td(\E^\vee)=(1-e^{-t})^d=t^d$ in the rational Chow ring $A(\PP^d)_\Q\cong\Q[t]/\ideal{t^{d+1}}$.
\end{proof}

\begin{proof}[Proof of Theorem~\ref{th:Chern}:]
Proposition~\ref{prp:genExact-proj} gives an exact sequence
\begin{equation}\label{50}
\xymat@C-0.6em{
0\ar[r] & \Omega^0(\PP\A)\ar[r]^-{\omega_\lambda} & \Omega^1(\PP\A)\ar[r]^-{\omega_\lambda} & \cdots\ar[r]^-{\omega_\lambda} & \Omega^{\ell-1}(\PP\A)\ar[r]^-{\omega_\lambda} & \F\ar[r] & 0
}
\end{equation}
for generic $\lambda$ with $|\lambda|=0$.
By \cite{OT95} and using \eqref{49}, 
\[
\F=(\O_{\PP V}/I_\lambda(\PP\A))(n-\ell)
\]
is a coherent sheaf supported on points of length equal to the {\em $\beta$-invariant} of $\A$,
\begin{equation}\label{45}
B=B(\A):=(-1)^{\ell-1}\pi(\PP\A,-1)
=\sum_{i=0}^{\ell-1}(-1)^{\ell-1-i}b_i.
\end{equation}
Using Lemma~\ref{lem:ch-pt}, then
\begin{equation}\label{53}
\ch(\F)=B t^{\ell-1}
\end{equation}
in $A(\PP V)_\Q$.

Recall from Propositions~\ref{prp:self-dual}.\eqref{sd3}, \ref{prp:proj}, and \ref{prp:proj-D}, that 
\[
\Omega^{\ell-2}(\PP\A)\cong\Omega^1(\PP\A)^\vee(n-\ell).
\]
Since $\PP\A$ is assumed to be locally tame, $\Omega^1(\PP\A)$ has a resolution of the form \eqref{eq:resM}.
Let $\alpha$ and $\beta$ be the Chern roots of $\F_0$ and $\F_1$, respectively.
Using the argument from the proof of Proposition~\ref{prp:N} (or from the proof of \cite[Prop.~2.6]{Har80}), and Lemma~\ref{lem:ch-pt}, we compute the Chern character
\begin{align}\label{52}
\ch_t(\Omega^{\ell-2}(\PP\A))&=(\ch_{-t}(\Omega^1(\PP\A))+Nt^{\ell-1})
e^{(n-\ell)t}\\
\nonumber&=\ch_{-t}(\Omega^1(\PP\A))e^{(n-\ell)t}+Nt^{\ell-1},
\end{align}
where we set $N:=N(\PP\A)$.
Similarly $\Omega^{\ell-1}(\PP\A)\cong\O_{\PP V}(n-\ell)$ gives
\begin{equation}\label{80}
\ch_t(\Omega^{\ell-1}(\PP\A))=e^{(n-\ell)t}.
\end{equation}
By restriction to a generic hyperplane and \cite[Rem.~3.2.3.(c)]{FultonBook} this gives
\begin{equation}\label{81}
c_1(\Omega^1(\PP\A))=n-\ell.
\end{equation}
Now the two equalities of Proposition~\ref{prp:chtop} give
\begin{align}
\label{82}L^{\ell-2}_{\alpha,\beta}&=\ch_t\Omega^{\ell-2}(\PP\A)-Nt^{\ell-1}, &&\text{by \eqref{52}, and}\\
\label{83} L^{\ell-1}_{\alpha,\beta}&=e^{c_1(\Omega^1(\PP\A))t}=
\ch_t(\Omega^{\ell-1}(\PP\A)),&&\text{by \eqref{81}, \eqref{80}.}
\end{align}

By Proposition~\ref{prp:Lebelt-proj} and \eqref{83},
\begin{equation}\label{51}
\Omega^p(\PP\A)=\bigwedge^p\Omega^1(\PP\A),\quad p\ne\ell-2.
\end{equation}
Now we can compute an Euler characteristic of \eqref{50}, using \eqref{51}. 
Proposition~\ref{prp:wedges} together with Theorem~\ref{thm:BorSer} yields
\begin{align}\label{54}
\ch_t(\F)
&=(-1)^{\ell-1}\sum_{p=0}^{\ell-1}(-1)^p\ch_t(\Omega^p(\PP\A))\\
\nonumber&=(-1)^{\ell-1}\sum_{p=0}^{\ell-1}(-1)^p
L^p_{\alpha,\beta}+L^{\ell-2}_{\alpha,\beta}-\ch_t(\Omega^{\ell-2}(\PP\A))\\
\nonumber&=c_{\ell-1}(\Omega^1(\PP\A))t^{\ell-1}-Nt^{\ell-1},
\text{ using \eqref{82}.}
\end{align}
Now take the coefficient of $t^{\ell-1}$ in \eqref{53} and \eqref{54} and use \cite[Exa.~3.2.2]{FultonBook} to find 
\[
B=\sum_{i=0}^{\ell-1}(-1)^{\ell-1-i}c_i-N.
\]
Finally, use Proposition~\ref{prp:slice-Chern} for $k=\ell-2$ to see $c_i=b_i$ for $0\leq i\leq \ell-2$.
It follows from the definition of $B$ in \eqref{45} that $c_{\ell-1}=b_{\ell-1}+N$, which completes the proof.
\end{proof}

\begin{exa}\label{ex:ziegler}
Let $\Z_1$ and $\Z_2$ be the two arrangements of $9$ lines in $\PP^2$ introduced in \cite[Ex.~8.7]{Zi89} and independently in \cite[Ex.~2.2]{Yu93}.
The two arrangements have isomorphic intersection lattices; however, as Schenck notes in \cite{Sc09}, they are distinguished by the property that the six triple points in $\Z_2$ lie on a conic in $\PP^2$, while the triple points in $\Z_1$ do not.  By formula~\ref{eq:locfreeChern}, we see $c_t(\Omega^1(\PP\Z_i)(1))=1+8t+22t^2$, the Poincar\'e polynomial, for $i=1,2$.
On the other hand, the Hilbert series of $\Omega^1_0(\A)$ is not combinatorially determined, for arrangements of rank at least $3$:
$h(\Omega^1_0(\Z_i),t)$ differs for $i=1,2$, as noted in \cite{Zi89}.  

Add a generic hyperplane to each arrangement above, to obtain two combinatorially equivalent arrangements $\Z_1^+$ and $\Z_2^+$ of $10$ planes in $\PP^3$.  
Then $\pi(\Z_i^+,t)=(1+8t+22t^2)(1+t)=1+9t+30t^2+22t^3$, and $\Z_i^+$ has a single non-free closed subarrangement of rank $3$, which is $\Z_i$. 
Computing with Macaulay~2~\cite{M2} shows $N(\Z_1)=20$, while $N(\Z_2)=22$. 
By Theorem~\ref{th:Chern}, we find
\[
c_t(\Omega^1(\PP\Z_i^+)(1))=
\begin{cases}
1+9t+30t^2+42t^3 & \text{for }i=1,\\
1+9t+30t^2+44t^3 & \text{for }i=2.
\end{cases}
\]
So we see that, for arrangements of rank at least $4$, the Chern polynomial of $\Omega^1(\PP\A)$ is not combinatorially determined, either.
\end{exa}


\begin{ack}
The authors would like to thank Hal Schenck, Uli Walther and Max Wakefield for helpful conversations, as well as the American Institute of Mathematics SQuaREs program for its hospitality and for facilitating this collaboration.
\end{ack}

\bibliographystyle{amsalpha}
\bibliography{aim}
\end{document}